\documentclass[conference]{IEEEtran}
\IEEEoverridecommandlockouts
\usepackage{cite}
\usepackage{amsmath,amssymb,amsfonts}
\usepackage{algorithm}
\usepackage{algpseudocode}
\usepackage{graphicx}
\usepackage{textcomp}
\usepackage{xcolor}
\usepackage[utf8x]{inputenc}
\usepackage{enumerate}
\usepackage{mathtools}
\usepackage{afterpage}
\DeclarePairedDelimiter\ceil{\lceil}{\rceil}

\def\BibTeX{{\rm B\kern-.05em{\sc i\kern-.025em b}\kern-.08em
    T\kern-.1667em\lower.7ex\hbox{E}\kern-.125emX}}
    
\algdef{SE}[SUBALG]{Indent}{EndIndent}{}{\algorithmicend\ }%
\algtext*{Indent}
\algtext*{EndIndent}

\begin{document}

\title{Hybrid Quantum Benders' Decomposition For Mixed-integer Linear Programming
}

\author{\IEEEauthorblockN{Zhongqi Zhao$^{1,2}$, Lei Fan$^2$, and Zhu Han$^{1}$}
\IEEEauthorblockA{\textit{Dept. of Electrical and Computer Engineering$^{1}$ and Dept. of Engineering Technology$^{2}$} \\
\textit{University of Houston, TX, USA}}}

\maketitle

\begin{abstract}
The Benders' decomposition algorithm is a technique in mathematical programming for complex mixed-integer linear programming (MILP) problems with a particular block structure. The strategy of Benders' decomposition can be described as a strategy of divide and conquer. The Benders' decomposition algorithm has been employed in a variety of applications such as communication, networking, and machine learning. However, the master problem in Benders' decomposition is still NP-hard, which motivates us to employ quantum computing. In the paper, we propose a hybrid quantum Benders' decomposition algorithm. We transfer the Benders' decomposition’s master problem into the quadratic unconstrained binary optimization (QUBO) model and solve it by the state-of-the-art quantum annealer. Then, we analyze the computational results and discuss the feasibility of the proposed algorithm. Due to our reformulation in the master problem in Benders' decomposition, our hybrid algorithm, which takes advantage of both classical and quantum computers, can guarantee the solution quality for solving MILP problems.
\end{abstract}


\section{Introduction}
In recent years, mixed-integer linear programming (MILP) is widely employed in many fields includes but is not limited to planning and scheduling \cite{canto2008application}, transportation and telecommunications \cite{costa2005survey}, energy and resource management \cite{cai2001solving}\cite{zhang2006hydro}, network design \cite{fortz2009improved}, and vehicle routing \cite{correa2007scheduling}. However, because MILP is an NP-hard problem \cite{bulut2021complexity} \cite{papadimitriou1982complexity}, it is desirable to have powerful solutions to solve large-scale MILP problems. 

Benders' decomposition algorithm \cite{geoffrion1972generalized} provides an efficient way to solve the MILP problem. In recent years, the Benders' decomposition method has become a popular algorithm. It exploits the structure of the MILP problem and successfully reduces the computation workload. Benders' decomposition divides a MILP problem into a master problem and a subproblem. The subproblem is a linear programming model that has strong duality. During the iterative solving process, the optimality cuts and feasibility cuts obtained from the solution of the subproblem will be added to the master problem. The iterative process will stop when the gap between the upper bound provided by the master problem and the lower bound provided by the subproblem is sufficiently small. However, the master problem in Benders' decomposition is still NP-hard, which motivates us to employ quantum computing. 

Due to the characteristics of parallel computing, the quantum computer can simultaneously calculate the final result with every possible model input. Therefore, quantum computers are often considered to have surpassed the computational speed of traditional computers. The quantum advantage has been proved by quantum algorithms, including Deutch-Jozsa Algorithm \cite{deutsch1992rapid}, Grover’s Algorithm \cite{grover1996fast}, Shor’s Algorithm \cite{shor1999polynomial}. In recent years, leading quantum computer companies such as D-Wave, IBM, Google, IonQ, Honeywell have made significant progress in hardware design. Normally, their quantum computers are in one of the two groups (i.e., analog quantum model, universal quantum gate model). D-Wave currently provides the quantum computer with the largest number of qubits on the market. 


By deploying the Ising model, the D-Wave's quantum annealer computer can solve the problem formulated by the quadratic unconstrained binary optimization (QUBO) model. It depicts the energy state with coupling qubits interaction and externally applied fields \cite{ising1925beitrag}. 

The power of quantum computers and challenges in MILP problems inspire us to design a hybrid quantum Benders' decomposition algorithm by jointly using quantum computing and classical computing techniques. However, there are several obstacles when we try to use Benders' composition to combine quantum and traditional computing. The first difficulty is how to convert the problem into an integer linear programming (ILP) problem recognized by the quantum computer. The Second difficulty is how to convert the NP-hard ILP problem into a QUBO model as an input to a quantum computer. In addition, the third difficulty is how to reformat the output of quantum computers from the binary solution to the decimal numeral system.

To overcome the above challenges, this paper reformulates the master problem as an ILP model in Benders' decomposition under the MILP problem. We transfer the master problem, which is an ILP model, into a QUBO form. Although the master problem in the Benders' decomposition is still NP-hard, we overcome this issue by proposing a hybrid Benders' decomposition algorithm with the D-Wave hybrid quantum computer for the MILP problem. Finally, we study a simple case to evaluate the performance of the D-Wave hybrid solver in solving the MILP problem. The contributions of this paper are summarized as follows.

\begin{itemize}
    \item We propose a hybrid quantum Benders' decomposition algorithm to find the solution for the MILP problem. Our hybrid quantum Benders' decomposition algorithm converges and returns the correct final result as the classical algorithm does.
    \item We propose an ILP model for the master problem in Benders' decomposition for the MILP problem. We reformulate the ILP model as a QUBO model recognized by the quantum annealing machine.
    \item We employ the quantum computer provided by D-Wave to solve the MILP problem by their Leap™ quantum cloud service. Our experiments show the possibility of using quantum computing to solve the MILP.
\end{itemize}

The rest of this paper is organized as follows. Section \ref{par:bdfm} introduces the basics about the MILP problem and the Benders’ decomposition. Section \ref{par:hqbd} illustrates our hybrid quantum Benders' Decomposition algorithm. Section \ref{par:NV} validates our algorithm by showing the corresponding simulation. Finally, Section \ref{par:C} concludes the whole paper.

\section{MILP and Benders' Decomposition Basics}
\label{par:bdfm}
\subsection{Mixed-integer Linear Programming}
MILP has been widely adopted optimization problems that include but are not limited to communication and networks. In the field of communication, for instance, it has been employed to the resource allocation in wavelength division multiple access (WDMA) \cite{elgamal2021q}, multi-access edge computing (MEC) \cite{yu2021uav}, offloading for fog computing \cite{wu2020energy}, etc. Furthermore, it plays an essential role in the field of network research, such as minimizing the network delay \cite{xuan2021minimizing}, scheduling virtual network re-configurations \cite{pan2021scheduling}, and improving the supply chain network \cite{santander2020closed}. Besides civilian applications, it is no doubt that MILP has been adopted for military purposes such as assisting military pilot training \cite{mak2021simultaneous}. These demonstrate that MILP is a powerful tool for either classical problems or future popular topics such as communication networks and resource allocation.
Consider a MILP model as follows:
\begin{equation} \label{eq:Original}
\begin{aligned}
\max_{x,y} \quad & c^{\intercal} x + h^{\intercal} y\\
\textrm{s.t.} \quad & Ax + Gy \leq b ,\\
  & x \in X,\ x \in \left\{0,1\right\}^{n} ,\ y \in \mathbb{R}^{p}_{+} .\\
\end{aligned}
\end{equation}
Here $x$ is a vector with binary variables and $A$ is its corresponding coefficient matrix in constraints. $y$ is a vector with non-negative continuous variables and $G$ is its corresponding coefficient matrix in constraints. $c^{\intercal}$ is the coefficient vector of variable $x$ in the objective function and $h^{\intercal}$ is the coefficient vector of variable $y$ in the objective function.

\begin{figure}[t]
\centering 
\includegraphics[width=8cm]{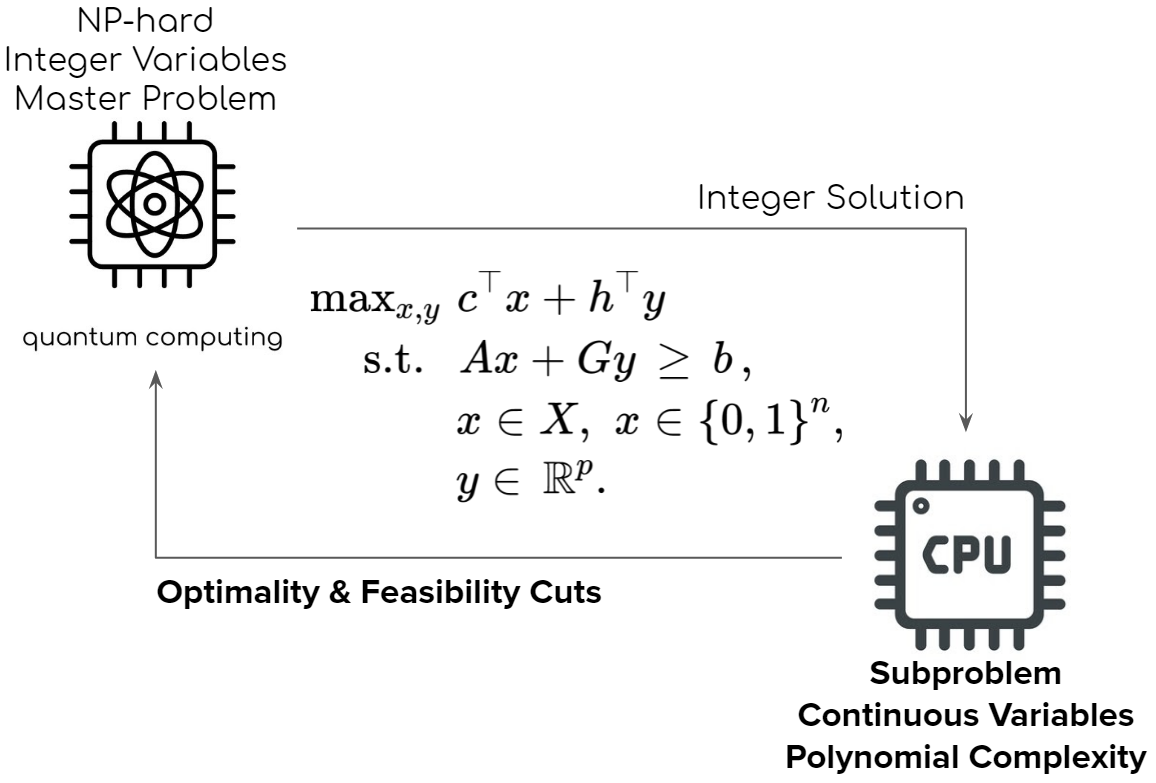}
\caption{The structure of hybrid Benders’ decomposition algorithm
for MILP using both quantum and traditional computing}
\label{fig:BD_structure}
\end{figure}

\subsection{Benders' Decomposition for MILP}
We first provide a brief introduction of the Benders' Decomposition. In the MILP problem (\ref{eq:Original}), for each possible choice of $\Bar{x} \in X$, we find the best choice for $y$ by solving a linear program. So we regard $y$ as a function of $x$. Then we replace the contribution of $y$ to the objective with a scalar variable representing the value of the best choice for a given $\Bar{x}$. We start with a crude approximation to the contribution of $y$ and then generate a sequence of dual solutions to tighten up the approximation. Hence, the original MILP problem in (\ref{eq:Original}) can be written equivalently to a master problem as follows. 

\begin{equation} \label{eq:relaxed Benders}
\begin{aligned}
\mbox{Master Problem:}\ \max_{x,t} \quad & c^{\intercal} x + t\\
\textrm{s.t.} \quad & \left(b-Ax\right)^{\intercal} u^{k} \geq t \quad \textrm{for} \ k \in \hat{K} ,\\
  & \left(b-Ax\right)^{\intercal} r^{j} \geq 0 \quad \textrm{for} \ j \in \hat{J}   ,\\
  & x \in X, \ x \in \left\{0,1\right\}^{n},\ t \in \mathbb{R} .\\
\end{aligned}
\end{equation}
We denote $K$ and $J$ as the extreme points $u^{k}$ and extreme rays $r^{k}$ of the dual polyhedron $Q = \left\{ u \in  \mathbb{R}^{m}_{+} \mid G^{\intercal}u \geq h \right\}$ generated by the linear programming’s duality, which is called the subproblem. The subproblem is as follows.
\begin{equation} \label{eq:Dual}
\begin{aligned}
\mbox{Subproblem:}\ z_{LP}(x) &= \min_{u} \ \left(b-Ax\right)^{\intercal}u\\
\textrm{s.t.} \quad & G^{\intercal}u \geq h ,\\
  & u \in \mathbb{R}^{m}_{+} .\\
\end{aligned}
\end{equation}
In the subproblem, if the inner product between $(b - Ax) $ and any dual ray $r^{j^{\prime}}$ is negative then $z_{LP}(x) = -\infty$. Equivalently, in this situation, the dual problem of Problem (\ref{eq:Dual}) is infeasible. Then, $x$ does not allow a feasible solution to the original mixed-integer problem (\ref{eq:Original}). Thus, we have a new feasibility cut, i.e.,
\begin{equation}
\label{eq:feasible}
    (b-Ax)^{\intercal}r^{j^{\prime}} \geq 0, \quad \ j^{\prime}\in J.
\end{equation}
If $x$ satisfies (\ref{eq:feasible}), then we yield an extreme point $u^{k^{\prime}}$ and the value of $z_{LP}(x)$ is given by
\begin{equation}
\label{eq:optimality}
    z_{LP}(x) = \min_{k^{\prime} \in K} \ \left(b-Ax\right)^{\intercal}u^{k^{\prime}}.
\end{equation}
Thus, problem (\ref{eq:Original}) is written equivalently as problem (\ref{eq:relaxed Benders}) by denoting $z_{LP}(x)$ as a continuous real number variable $t$. Equation (\ref{eq:optimality}) posts a new optimality cut in master problem. In addition to the extreme points and extreme rays, we use $\hat{K}$ and $\hat{J}$ to denote the current known extreme points and extreme rays of $Q$, respectively, where $\hat{K} \subseteq K$ and $\hat{J} \subseteq J$.

\section{Hybrid Quantum Benders' Decomposition}
\label{par:hqbd}
Quantum computer is a powerful tool to solve NP-hard integer problems. Accordingly, we introduce the quantum computer to solve the master problem in (\ref{eq:relaxed Benders}) since it has a special formulation that can convert to an ILP problem. The subproblem is a linear programming model which can be solved well on the classical computer. As Fig. \ref{fig:BD_structure} shows, 
we design a hybrid algorithm in which the master problem is solved by the quantum computing techniques and the subproblem is solved by the classical computing techniques. The two different computers are communicating with each other. The subproblem returns new optimality and feasibility cuts to tighten the bounds of the master problem. The integer solution returned from the master problem gives the subproblem a direction to find new optimality and feasibility cuts.

\begin{table}[t]
\renewcommand\arraystretch{2}
\caption{Table of Common Constraint-Penalty Pairs}
\noindent \vspace{-0.4cm}
\label{table2}
\begin{center}

\begin{tabular}{|l|l|}
\hline
\textbf{Constraint}& \textbf{Equivalent Penalty}\\
\hline
\label{table:penalty form}
$x_{1} + x_{2} = 1 $& $P(x_{1} + x_{2} - 1)^{2} $ \\
\hline
$x_{1} + x_{2} \geq 1 $& $P(1 - x_{1} - x_{2} + x_{1}x_{2} )^{2} $\\
\hline
$x_{1} + x_{2}  \leq 1 $& $P(x_{1}x_{2}) $ \\
\hline
$x_{1} + x_{2} + x_{3}\leq 1 $& $P(x_{1}x_{2} + x_{1}x_{3} + x_{2}x_{3}) $\\
\hline

\end{tabular}
\end{center}
\end{table}

\subsection{Quantum Formulation}
Quantum annealers are able to solve the optimization problem in a QUBO formulation. To leverage state-of-art quantum annealers provided by D-Wave, the ILP problem has to be converted to the corresponding QUBO formulation. The definition of QUBO are (\ref{eq:QUBO}) and (\ref{eq:QUBO2}). Let $f$:$ \left\{0,1\right\}^{n} \rightarrow \mathbb {R}$ be a quadratic polynomial over binary variables $x$ of length $n$ and $q_{ij}$ is corresponding cost coefficient of $x_{i}x_{j}$. We get equation (\ref{eq:QUBO})

\begin{equation} \label{eq:QUBO}
f_{Q}\left(\mathbf{x}\right)= \sum_{i=0}^{n-1}\sum_{j=0}^{n-1}q_{ij}x_{i}x_{j} = \mathbf{x}^{\intercal} \mathbf{Q} \mathbf{x}.
\end{equation}

The QUBO problem consists of finding a binary vector $\mathbf{x}^{*}$ that is minimal with respect to $f$ among all other binary vectors, namely,
\begin{equation}\label{eq:QUBO2}
    \mathbf{x}^{*} = {\arg\min}_{\mathbf{x}} \ f\left(\mathbf{x}\right).
\end{equation}
Here $\mathbf{x}$ is a vector of binary variables where the length of $n$, and $\mathbf{Q}$ is either an upper-diagonal matrix or a symmetric matrix. Since (\ref{eq:QUBO}) is an unconstrained optimization model, we need to reformulate our constrained ILP as unconstrained QUBO by using penalties. Next, we get the optimal solution by finding the best penalty coefficients of the constraints. The principles of transforming classical constraints to their equivalent penalties are displayed in the TABLE \ref{table:penalty form}. Here $x_1$, $x_2$ and $x_3$ are binary variables.
$s_i$ is a binary slack variable. $a_l$ is the coefficient for the corresponding slack variable. $b$ is a constant. $P$ is a user-defined penalty coefficient.

\subsection{Variable Representation}
Now consider problem (\ref{eq:relaxed Benders}), the initial binary decision variables make up the vector $\mathbf{x} \in X$ with length of $n$. $\Bar{t} \in \mathbb{R} $. In order to reformulate the master problem into the QUBO formulation, we needs to represent the continuous variable using binary bits. We use a binary vector $\mathbf{w}$ with length of $M$ bits to replace continuous variable $\Bar{t}$. In general, $\Bar{t}$ requires the binary numeral system assigning $M$ bits to replace continuous variable $\Bar{t}$. 
Then we can recover the $\Bar{t}$ by 

\begin{equation} \label{eq: finaltrepresentation}
    \begin{aligned}
        \Bar{t}  &=  \sum^{\Bar{m}_{+}}_{i = - \underline{m}} 2^{i}w_{\left(i+\underline{m}\right)} - \sum^{\Bar{m}_{-}}_{j = 0} 2^{j}w_{j+ 1 + \underline{m} + \Bar{m}_{+} } \\
        &= \Bar{t}\left(\mathbf{w}\right). \\
    \end{aligned}
\end{equation}
In (\ref{eq: finaltrepresentation}), $\Bar{m}_{+}$ is the number of bits that assigned to represent positive integer part, $\underline{m}$ is the number of bits assigned to represent the positive decimal part, and $ \Bar{m}_{-} $ is the number of bits that represent the negative integer part.

\subsection{QUBO Setup}
Let's reconsider problem (\ref{eq:relaxed Benders}) and replace $\Bar{t}$ in (\ref{eq: finaltrepresentation}). The new problem is only depend on binary decision variables $\mathbf{x}$ and $\mathbf{w}$. Then, the new master problem is

\begin{equation} 
    \begin{aligned} 
        \max_{\mathbf{x},\mathbf{w}} \quad & c^{\intercal} \mathbf{x} + \sum^{\Bar{m}_{+}}_{i = -\underline{m}} 2^{i}w_{\left(i+\underline{m} \right)} - \sum^{\Bar{m}_{-}}_{j = 0} 2^{j}w_{j+ 1 + \underline{m} + \Bar{m}_{+} }\\
        \textrm{s.t.} \quad &  \left(b-A\mathbf{x}\right)^{\intercal} u^{k} \geq \Bar{t}\left(\mathbf{w}\right), \quad \textrm{for} \ k \in \hat{K},\\
          &  \left(b-A\mathbf{x}\right)^{\intercal} r^{j} \geq 0, \quad \textrm{for} \ j \in \hat{J},\\
          & \mathbf{x} \in X ,\quad \mathbf{x} \in \left\{0,1\right\}^{n},   \\
          & \mathbf{w} \in W ,\quad \mathbf{w} \in \left\{0,1\right\}^{M} .  \\
    \end{aligned}\label{eq:trasnferred s}
\end{equation}

Then we create a new binary decision variable collection $\mathbf{x}^{\prime} = \left\{\mathbf{w},\mathbf{x}\right\}$ to set up our QUBO matrix, where the length of the binary vector $\mathbf{x}^{\prime}$ is $n + M$.
According to the rule of setting up QUBO, we will transfer our objective and constraints one by one to be a symmetric matrix.
\subsubsection{Objective Function}
    Because the quantum computer only accepts a quadratic polynomial over binary variables. Following the principle of the QUBO formulation, we convert the objective function to be a QUBO matrix as follows. 
    \begin{equation} 
    \begin{aligned} \label{eq:QUBO_Obj}
    &c^{\intercal} \mathbf{x} + \sum^{\Bar{m}_{+}}_{i = -\underline{m}} 2^{i}w_{\left(i+\underline{m} \right)} - \sum^{\Bar{m}_{-}}_{j = 0} 2^{j}w_{j+ 1 + \underline{m} + \Bar{m}_{+} } \\ 
    \Rightarrow & \mathbf{Q}_{obj} = \sum^{\Bar{m}_{+}}_{i = -\underline{m}}w_{\left(i+\underline{m} \right)}2^{i} w_{\left(i+\underline{m} \right)} + \mathbf{x}^{\intercal}\textrm{diag}\left( c \right)\mathbf{x} \\
    & - \sum^{\Bar{m}_{-}}_{j = 0} w_{j+ 1 + \underline{m} + \Bar{m}_{+} } 2^{j}w_{j+ 1 + \underline{m} + \Bar{m}_{+} }.
    \end{aligned}
\end{equation}

\subsubsection{Optimality Cuts}
    
Similarly, following the principle of constraint-penalty pairs in TABLE \ref{table:penalty form}, we not only introduce the penalty but also convert the optimality cuts constraint to be a QUBO matrix as follow. 
\begin{equation*} 
\begin{aligned} \label{eq:QUBO_C1}
&\Bar{t}\left(\mathbf{w}\right) + \left(u^{k}\right)^{\intercal}A\mathbf{x}  \leq b^{\intercal} u^{k}, \ \textrm{for} \ k \in \hat{K}.\\
    \Rightarrow & P_{k}\left( \Bar{t}\left(\mathbf{w}\right) + \left(u^{k}\right)^{\intercal}A\mathbf{x} + \sum_{l = 0}^{\Bar{l}^{K}} 2^{l} s^{K}_{kl} - b^{\intercal} u^{k} \right)^{2},\\
        \textrm{where} \   &\bar{l}^{K} = \ceil*{\log_{2}\left(b^{\intercal} u^{k} - \min_{\mathbf{w},\mathbf{x}}\left( \Bar{t}\left(\mathbf{w}\right) + \left(u^{k}\right)^{\intercal}A\mathbf{x}\right)\right)}.
    \end{aligned}
\end{equation*}
    
\subsubsection{Feasibility Cuts}
    Similarly, following the principle of constraint-penalty pairs in TABLE \ref{table:penalty form}, we convert the feasibility cuts constraint to be a QUBO matrix with penalties as follows:
\begin{equation*} 
    \begin{aligned} \label{eq:QUBO_C2}
    &\left(r^{j}\right)^{\intercal}A\mathbf{x} \leq b^{\intercal} r^{j}, \quad \textrm{for} \ j \in \hat{J}.\\
        \Rightarrow & P_{j}\left( \left(r^{j}\right)^{\intercal}A\mathbf{x} + \sum_{l = 0}^{\Bar{l}^{J}} 2^{l} s^{J}_{kl} - b^{\intercal} r^{j} \right)^{2},\\
        \textrm{where} \   &\bar{l}^{J} = \ceil*{\log_{2}\left(b^{\intercal} r^{j} - \min_{\mathbf{x}}\left(  \left(r^{j}\right)^{\intercal}A\mathbf{x}\right)\right)}.
    \end{aligned}
\end{equation*}
\par    

\subsection{Proposed Algorithm}
We deploy the D-Wave solver in our proposed heuristic algorithm to solve the model. In addition, we need to carefully tune the penalties for a decent QUBO model. An extremely large penalty may cause quantum annealer malfunctioning since it will explode the coefficients. Similarly, a soft penalty may make quantum annealer ignore the corresponding constraints. As long as the penalty is well-tuned, the quantum solver will give the right answer with a relatively high probability. Therefore, we summarize our proposed hybrid quantum Benders' decomposition algorithm as \textbf{Algorithm \ref{alg:Q_bender}}.
\begin{algorithm} 
\caption{Hybrid Quantum  Benders' Decomposition Algorithm}
\label{alg:Q_bender}
\begin{algorithmic}
\Require Initial sets $\hat{K}$ of extreme points and $\hat{J}$ of extreme rays of $Q$
\State $\Bar{t}$ $\gets$ $+\infty$ 
\State $\underline{t}$ $\gets$ $-\infty$ 

\While{$ \mid \Bar{t} - \underline{t} \mid \geq \epsilon$}

\State $\mathbf{P}$ $\gets$ Appropriate penalties numbers or arrays
\State $\mathbf{Q}$ $\gets$ Reformulate both objective and constraints in ($\ref{eq:relaxed Benders}$) and construct the QUBO formulation by using corresponding rules
\State $\mathbf{x}^{\prime}$ $\gets$ Solve problem (\ref{eq:QUBO}) by quantum computer.
\State $\Bar{t}$ $\gets$ Extract $\mathbf{w}$ and replace the $\bar{t}$ with $\Bar{t}\left(\mathbf{w}\right)$ (\ref{eq: finaltrepresentation})
\State $z_{LP} (\mathbf{x})$ $\gets$ Solve the problem (\ref{eq:Dual})
\State $\underline{t}$ $\gets$ $z_{LP} (\mathbf{x})$
\If{$ z_{LP} (\mathbf{x}) = -\infty$ }
    \State An extreme ray $j$ of $Q$ has been found.
    \State $\hat{J}$ = $\hat{J}$ $\cup$ $\left\{j\right\}$
    
\ElsIf{$z_{LP} (\mathbf{x}) < \Bar{t}$ \textbf{and} $\Bar{t} \neq +\infty$ }
    \State An extreme point $k$ of $Q$ has been found.
    \State $\hat{K}$ = $\hat{K}$ $\cup$ $\left\{k\right\}$
\EndIf
\EndWhile\\
\Return $\Bar{t}$, $\mathbf{x}$
\end{algorithmic}
\end{algorithm}

\section{NUMERICAL VALIDATION}
\label{par:NV}

We validate the proposed algorithm by running on a hybrid D-Wave quantum processing unit (QPU). Because the input of interest to practice is too large to fit onto current-model QPUs and be solved directly by quantum annealing. We choose to use the hybrid model instead of the pure QPU. The reason for using the hybrid model is that the hybrid solver overcomes a lot of input size barriers and allows the QPU to accept a large input. More details can be found in \cite{dwave}. We accessed the D-Wave system by Leap™ quantum cloud service.

\subsection{Example Setup}
In our simulation, we consider a simple problem to test our proposed quantum algorithm, where $\mathbf{x}\in X$, $\mathbf{x} \in \left\{0,1\right\}^{2} $, $\mathbf{y} \in Y$, $\mathbf{y} \geq 0$

\begin{equation*}
A = 
\begin{bmatrix}
0 & 0 \\
0 & 0 \\
0 & 0 \\
0 & 0 \\
-1 & -1 \\
-1 & 0 \\
-1 & 0 \\
0 & -1 \\
0 & -1 \\
\end{bmatrix}
,\
G = 
\begin{bmatrix}
1 & 0 & 1 & 0 \\
1 & 0 & 0 & 1 \\
0 & 1 & 1 & 0 \\
0 & 1 & 0 & 1 \\
0 & 0 & 0 & 0 \\
1 & 0 & 0 & 0 \\
0 & 1 & 0 & 0 \\
0 & 0 & 1 & 0 \\
0 & 0 & 0 & 1 \\
\end{bmatrix}
,\
b = 
\begin{bmatrix}
1 \\
1 \\
1 \\
1 \\
-1 \\
0 \\
0 \\
0 \\
0 \\
\end{bmatrix},
\end{equation*}
\begin{equation*}
h^{\intercal} = 
\begin{bmatrix}
8 & 9 & 5 & 6\\
\end{bmatrix},\quad
c^{\intercal} = 
\begin{bmatrix}
-15 & -10\\
\end{bmatrix}.
\end{equation*}

\subsection{Result}

\begin{figure}[t]
\centering 
\includegraphics[width=8.5cm]{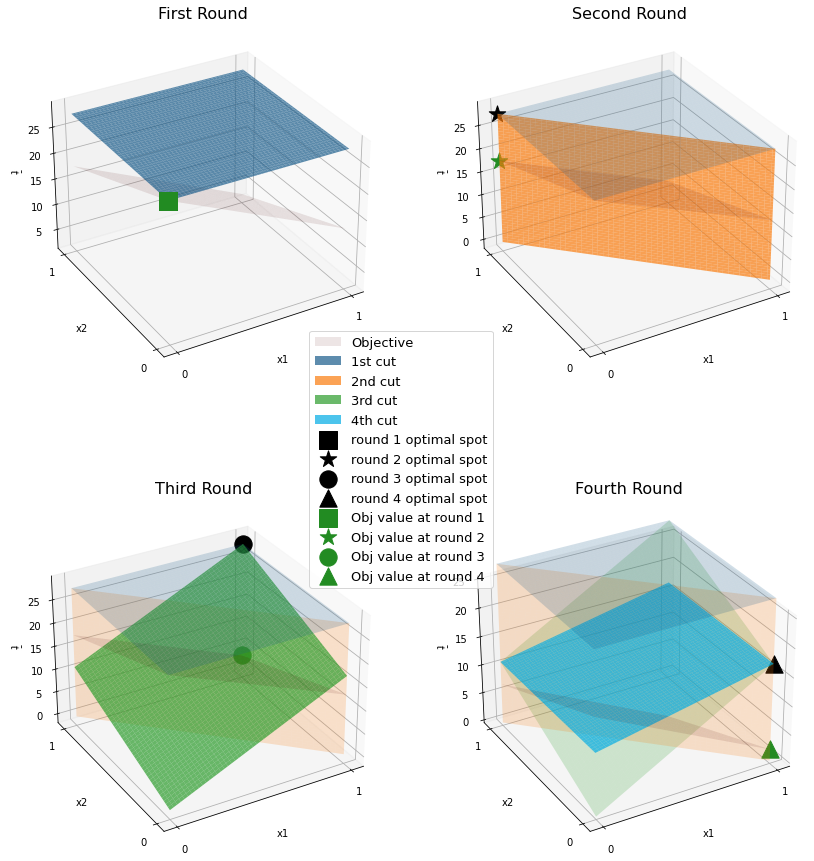}
\caption{Benders' Decomposition Process}
\label{fig:cutting}
\vspace*{-4mm}
\end{figure}

Fig. \ref{fig:cutting} shows how our proposed hybrid quantum Benders' decomposition algorithm adds optimality and feasibility cuts to get the solution and corresponding $\Bar{t}$. In the example, the algorithm takes 4 rounds to let $t$ converge. Four cuts are added iteratively in the space and tighten the feasible region. The algorithm can find the optimal solution in each round. Therefore, Fig. \ref{fig:cutting} demonstrates that our algorithm is reliable and efficient. In Fig. \ref{fig:tutl}, the hidden dashed line denotes that the lower bound in the corresponding round is negative infinity. As we can see in the graph, the upper bound and lower bound converge. Our algorithm only takes two rounds to find the non-negative infinity lower bound. This result proves that our proposed algorithm is mathematically consistent with the classic Benders' decomposition algorithm. In other words, as long as the classic Benders' decomposition algorithm can solve it, our algorithm can at least achieve the same. 

\section{Conclusion}
\label{par:C}
In this paper, we not only reformulate the master problem of Benders' decomposition for the MILP problem as an ILP model but also successfully convert it to a QUBO model. our hybrid quantum Benders' decomposition algorithm converges and returns the correct final result as the classical algorithm does. In addition to that, our algorithm also guarantees the solution quality for solving the MILP problem. In our simulation, we solved an NP-hard MILP problem by using the hybrid quantum computer provided by D-Wave. Therefore, we can conclude that the quantum computer can potentially replace the classical computer in solving the NP-hard master problem of Benders' decomposition algorithm in MILP problems.

\afterpage{%
\begin{figure}[t]
\centering 
\includegraphics[width=8cm]{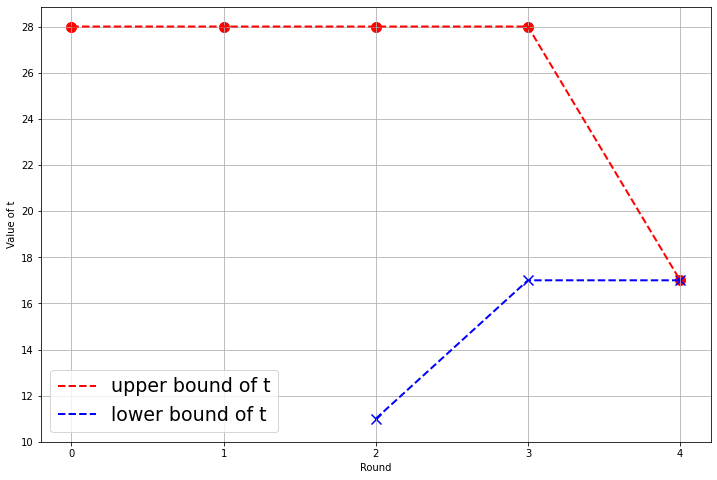}
\caption{Convergence of $\Bar{t}$ and $\underline{t}$}
\label{fig:tutl}
\vspace*{-2mm}
\end{figure}%
}

\bibliographystyle{IEEEtran}

\bibliography{ref_hybrid}

\end{document}